\documentclass{amsart}
\usepackage{ifthen, citesort}
\usepackage{amssymb,amsmath}
\usepackage{epic}
\input epsf
\usepackage{epsfig}

\DeclareMathOperator{\dist}{dist}
\def\R{\mathbb{R}}
\def\N{\mathbb{N}}

\def\fracp#1#2{\frac{\partial #1}{\partial #2}}
\def\theta{\vartheta}
\def\phi{\varphi}
\def\epsilon{\varepsilon}

\def\ctw{\cite{TrudWangAffineComplete}}

\newtheorem{theorem}{Theorem}

\theoremstyle{definition}

\newtheorem{defrem}[theorem]{Definition and Remark}

\theoremstyle{remark}
\newtheorem{remark}[theorem]{Remark}

\newcommand{\abs}[1]{\left\lvert#1\right\rvert}

\begin{document}

\title{Convex Functions with Unbounded Gradient}

\author{Oliver C. Schn\"urer}
\address{Oliver Schn\"urer, FU Berlin, Arnimallee 2-6, 
  14195 Berlin, Germany}
\curraddr{}
\email{Oliver.Schnuerer\fuhome}
\def\fuhome{@math.fu-berlin.de}
\thanks{The author is a member of SFB 647/3B ``Raum -- Zeit -- Materie''}


\date{April 2005, revised September 2005.}

\dedicatory{}


\begin{abstract}
We show that domains, that allow for convex functions with
unbounded gradient at their boundary, are convex.
\end{abstract}

\maketitle


In this paper, we prove the following
\begin{theorem}\label{einl thm}
Let $\Omega\subset\R^n$ be any domain. If there exists a 
convex function $u:\Omega\to\R$, such that $\nabla u$
becomes unbounded near $\partial\Omega$, then $\Omega$
is convex.
\end{theorem}

This problem arises in the context of affine hypersurfaces
\ctw, where a similar statement is proven.\\[0.5ex]
\indent{\bf Acknowledgement:} 
We want to thank Neil Trudinger, who brought this problem to
our attention and encouraged us to provide a geometric 
proof. We also want to thank Klaus Ecker and Free University
Berlin for support during the preparation of this paper
and Ann Bj\"orner for converting our sketch into a
beautiful illustration.\\[0.7ex]
\indent We will use the following notions of convexity.
\begin{defrem}
Let $\Omega\subset\R^n$ be a given domain.
\par 
A function $u:\Omega\to\R$ is called convex, if 
$$u(\tau x+(1-\tau)y)\le\tau u(x)+(1-\tau)u(y)$$
for all $x,\,y\in\Omega$ and all $\tau\in(0,\,1)$ such that
$\tau x+(1-\tau)y\in\Omega$.
\par
A convex function is locally Lipschitz continuous and thus
differentiable almost everywhere. We obtain
$$u(y)\ge u(x)+\langle\nabla u(x),\,y-x\rangle$$
for all $y\in\Omega$ and all $x\in\Omega$, where $u$ is
differentiable.
\par
A function $u$ is called locally convex, if
$$u(\tau x+(1-\tau)y)\le\tau u(x)+(1-\tau)u(y)$$
for all $x,\,y\in\Omega$ and $\tau\in(0,\,1)$ such that
$tx+(1-t)y\in\Omega$ for all $t\in(0,\,1)$.
\par
For $u\in C^2(\Omega)$, local convexity is equivalent to
positive semi-definiteness of its Hessian in $\Omega$.
\par
Note especially that locally convex functions do not need
to be convex, unless their domain of definition is a 
convex set, i.\,e.\ $x,\,y\in\Omega$, $\tau\in(0,\,1)$
$\Longrightarrow$ $\tau x+(1-\tau)y\in\Omega$.
\end{defrem}

In order to strengthen Theorem \ref{einl thm},
we fix $x_0\in\Omega$ and define $\Omega_1=\Omega_1(x_0)$ by
$$\Omega_1:=\{x\in\Omega:tx+(1-t)x_0\in\Omega
\text{~for~}0\le t\le1\}$$
as the largest subdomain of $\Omega$ that is star shaped 
with respect to $x_0$. Note that $\Omega$ is convex if and
only if $\Omega=\Omega_1(x)$ for all $x\in\Omega$.

For convenience, we will assume from now on that $x_0$ is the
origin. So the domain $\Omega_1=\Omega_1(0)$ is star shaped
with respect to the origin. 

We define
$$\partial_2\Omega_1:=\{x\in\partial\Omega_1:\exists\epsilon>0:
  tx\in\partial\Omega_1\forall t\in(1-\epsilon,\,1+\epsilon)\}$$
and $\partial_1\Omega_1:=
\partial\Omega_1\setminus\partial_2\Omega_1.$

For $\Omega$ as in Theorem \ref{einl thm}, we deduce that 
$\nabla u$ becomes unbounded near $\partial_1\Omega_1$.

Thus the following result generalizes Theorem \ref{einl thm}
and \cite[Lemma 2.4]{TrudWangAffineComplete}. 

\begin{theorem}\label{main thm}
Let $\Omega\subset\R^n$ be open and star shaped 
with respect to $0\in\Omega$. If there
exists a convex function $u:\Omega\to\R$, such that
$$\Vert\nabla u\Vert_{L^\infty(B_\delta(x)\cap\Omega)}=\infty$$
for any $x\in\partial_1\Omega$ and any $\delta>0$, 
then $\Omega$ is convex.
\end{theorem}

Here, the set $\partial_1\Omega$ is defined as above. We consider
$\nabla u$ only at those points, where $u$ is differentiable.

\begin{proof}
We argue by contradiction. Assume that there exist
$\tilde y_0$ and $\tilde y_1$ in $\Omega$ and $\tilde\tau\in(0,\,1)$
such that $$\tilde\tau\tilde y_0+(1-\tilde\tau)\tilde y_1\notin\Omega.$$
Let $\lambda>0$ be the supremum over all positive numbers such that
$$t\lambda\tilde y_0+(1-t)\lambda\tilde y_1\in\Omega\quad
\text{for~all~}0\le t\le1.$$
As $0$ is an interior point of $\Omega$, 
$\lambda$ is positive. Define $y_0:=\lambda\tilde y_0$ and
$y_1:=\lambda\tilde y_1$. Let $\tau\in(0,\,1)$ be such that
$\tau y_0+(1-\tau)y_1\in\partial\Omega$. According to the 
definition of $\lambda$,
$$\mu(\tau y_0+(1-\tau)y_1)\in\Omega\quad\text{for~all~}
0\le\mu<1.$$
So we deduce that $\tau y_0+(1-\tau)y_1\in\partial_1\Omega$.
Thus, there exist $x_i\in\Omega$, $i\in\N$, such that
$x_i\to\tau y_0+(1-\tau)y_1$ for $i\to\infty$, $\nabla u$
exists at $x_i$ for all $i$, and $\abs{\nabla u(x_i)}\to\infty$
as $i\to\infty$. This contradicts the gradient bounds that we
will prove in the following.
\par
As $\Omega$ is open and star shaped, there exists $\epsilon>0$
such that $B_{3\epsilon}(y_k)\subset\Omega$, $k=0,\,1$,
and $B_{3\epsilon}(0)\subset\Omega$.
Assume that the balls 
$B_{3\epsilon}(y_0)$, $B_{3\epsilon}(y_1)$,
$B_{3\epsilon}(\tau y_0+(1-\tau)y_1)$, and 
$B_{3\epsilon}(0)$ are disjoint. 
By scaling $u$, we may therefore arrange that $\abs u\le1$
in $B_{2\epsilon}(y_0)\cup B_{2\epsilon}(y_1)\cup
B_{2\epsilon}(0)$.
\par
Assume furthermore that, after an appropriate rotation, $y_1-y_0$
is a positive multiple of $e_1=(1,\,0,\,\ldots,\,0)$.
\begin{figure}[htb]
\epsfig{file=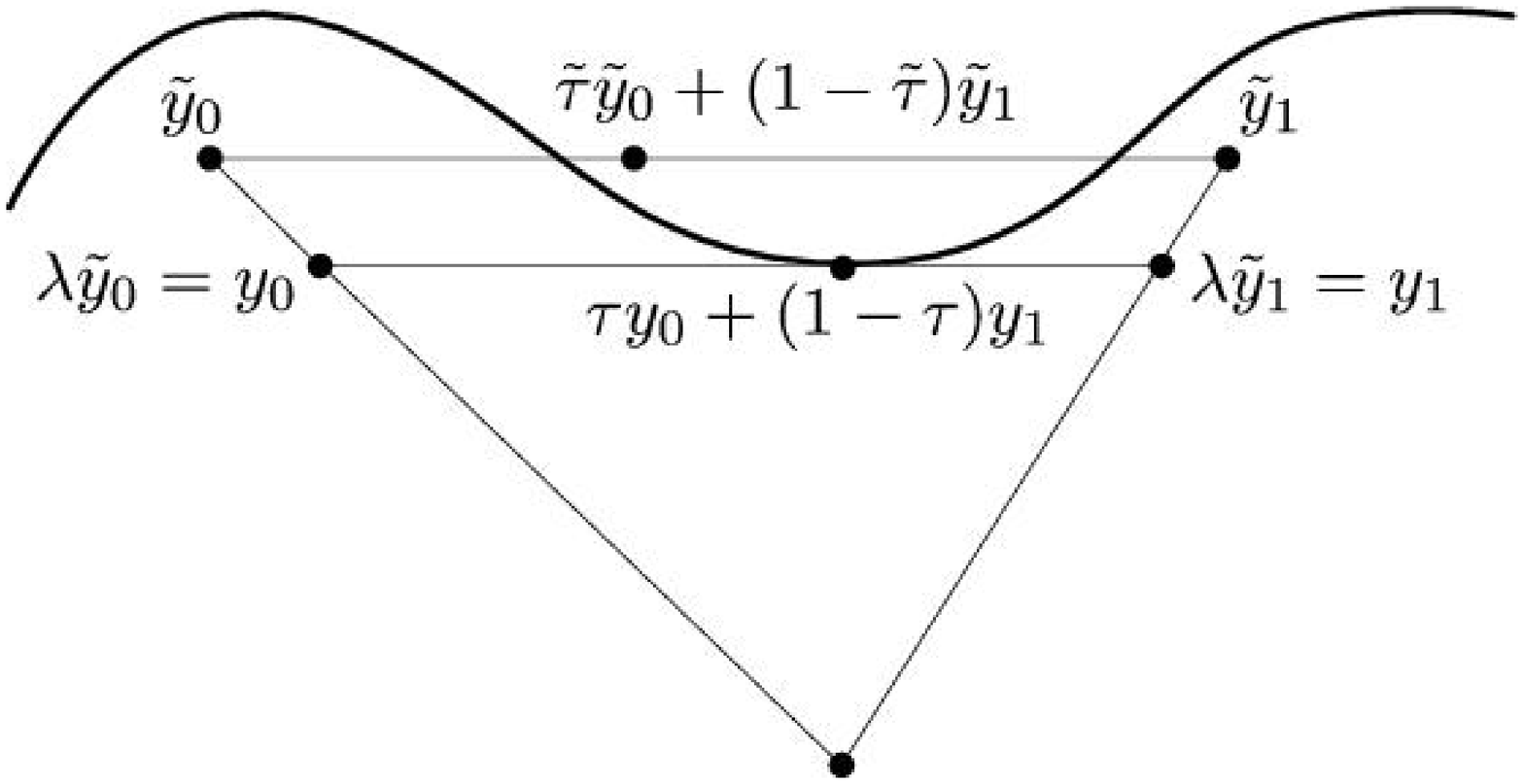, 
  width=0.9\textwidth}
\caption{Geometric situation} 
\label{notation}
\end{figure}
\par
We will now use the convexity of $u$ and boundedness of $u$ in balls
around $y_0$ and $y_1$ to prove explicit bounds on 
$\nabla u$ near $\tau y_0+(1-\tau)y_1$.
Set $u_k:=\fracp u{x^k}$, $1\le k\le n$. 
\par
The first step is to bound $u$ in
$B_\epsilon(\tau y_0+(1-\tau)y_1)$.
Let $z\in B_\epsilon(0)$ be such that
$\tau y_0+(1-\tau)y_1+z\in\Omega$. Convexity of $u$ and $u\le 1$
in $B_\epsilon(y_k)$, $k=0,\,1$, imply that
\begin{align*}
u(\tau y_0+(1-\tau)y_1+z)=&u(\tau(y_0+z)+(1-\tau)(y_1+z))\\
\le&\tau u(y_0+z)+(1-\tau)u(y_1+z)\\
\le&\tau\cdot 1+(1-\tau)\cdot 1=1.
\end{align*}
For a lower bound on $u$, we compare $u$ at 
$\tau y_0+(1-\tau)y_1+z=:p$, at $2\epsilon\frac p{\abs p}$, 
and at the origin. We get $\abs p\le\abs{\tau y_0+(1-\tau)y_1}
+\epsilon$,
$$2\epsilon\frac p{\abs p}=tp+(1-t)0$$
with $t=\frac{2\epsilon}{\abs p}$, and
$$\abs p\ge\abs{\tau y_0+(1-\tau)y_1}-\abs z\ge6\epsilon-\epsilon,$$
as $B_{3\epsilon}(\tau y_0+(1-\tau)y_1)\cap B_{3\epsilon}(0)=\emptyset$.
Thus 
$$u\left(2\epsilon\frac p{\abs p}\right)\le tu(p)+(1-t)u(0)$$
and we obtain that 
\begin{align*}
u(p)\ge&\frac1t\left(u\left(2\epsilon\frac p{\abs p}\right)
-(1-t)u(0)\right)\\
\ge&\frac{\abs{\tau y_0+(1-\tau)y_1}+\epsilon}{2\epsilon}(-1-1)\\
=&-\frac1\epsilon(\abs{\tau y_0+(1-\tau)y_1}+\epsilon).
\end{align*}
Thus $\abs u\le 1+\frac1\epsilon\abs{\tau y_0+(1-\tau)y_1}\equiv C_0$ 
in $B_\epsilon(\tau y_0+(1-\tau)y_1)$.
\par
The next step is to bound $\abs{u_1}$ in
$B_\epsilon(\tau y_0+(1-\tau)y_1)$. 
Remember that $y_1-y_0$ is a positive multiple of $e_1$.
As balls of radius $3\epsilon$
as chosen above are disjoint, we get
\begin{align*}
\abs{(1-\tau)(y_0-y_1)}=\abs{y_0-(\tau y_0+(1-\tau)y_1)}\ge&6\epsilon
\intertext{and}
\abs{\tau(y_0-y_1)}=\abs{(\tau y_0+(1-\tau)y_1)-y_1}\ge&6\epsilon.
\end{align*}
Thus we obtain
\begin{equation}\label{tau bound}
1-\tau\ge\frac{6\epsilon}{\abs{y_0-y_1}}\quad\text{and}\quad
\tau\ge\frac{6\epsilon}{\abs{y_0-y_1}}.
\end{equation}
Let $z\in B_\epsilon(0)$ and $\tau y_0+(1-\tau)y_1+z=:p$ be
an arbitrary point in $B_\epsilon(\tau y_0+(1-\tau)y_1)$,
where $u$ is differentiable.
We deduce that 
\begin{align*}
p+(1-\tau)(y_0-y_1)\in&B_\epsilon(y_0)
\intertext{and}
p+\tau(y_1-y_0)\in&B_\epsilon(y_1).
\end{align*}
Convexity and boundedness of u in $B_\epsilon(y_0)\cup
B_\epsilon(y_1)\cup B_\epsilon(\tau y_0+(1-\tau)y_1)$ imply that
\begin{align*}
\tau\abs{y_1-y_0}\langle\nabla u(p),\,e_1\rangle=&
\langle\nabla u(p),\,\tau(y_1-y_0)\rangle\\
\le&u(p+\tau(y_1-y_0))-u(p)\\
\le&1+C_0
\intertext{and}
(1-\tau)\abs{y_1-y_0}\langle\nabla u(p),\,-e_1\rangle=&
\langle\nabla u(p),\,(1-\tau)(y_0-y_1)\rangle\\
\le&u(p+(1-\tau)(y_0-y_1))-u(p)\\
\le&1+C_0.
\end{align*}
We insert the bounds \eqref{tau bound} and deduce that
$$\abs{\langle\nabla u(p),\,e_1\rangle}=\abs{u_1}(p)\le
\frac1{6\epsilon}(1+C_0).$$
\par
The last step is to bound $u_k$ in $B_\epsilon(\tau y_0+(1-\tau)y_1)$
for $2\le k\le n$. Consider an arbitrary vector $\xi\in\R^n$ such that
$\abs\xi=1$ and $\langle\xi,\,e_1\rangle=0$. As above, let
$p\in B_\epsilon(\tau y_0+(1-\tau)y_1)\cap\Omega$ be chosen
arbitrarily such that $u$ is differentiable in $p$. Then we have
$$p+(1-\tau)(y_0-y_1)\pm\epsilon\xi\in B_{2\epsilon}(y_0).$$
We remark that considering $p+\tau(y_1-y_0)\pm\epsilon\xi
\in B_{2\epsilon}(y_1)$ in the following yields a similar 
estimate. As $u$ is convex, we get
$$\langle\nabla u(p),\,(1-\tau)(y_0-y_1)\pm\epsilon\xi\rangle
\le u(p+(1-\tau)(y_0-y_1)\pm\epsilon\xi)-u(p).$$
So we deduce
\begin{align*}
\epsilon\langle\nabla u(p),\,\pm\xi\rangle\le&
u(p+(1-\tau)(y_0-y_1)\pm\epsilon\xi)-u(p)
+\langle\nabla u(p),\,(1-\tau)(y_1-y_0)\rangle\\
\le&\left(1+\frac{\abs{y_0-y_1}}{6\epsilon}\right)(1+C_0).
\end{align*}
Thus $\abs{\nabla u}$ is bounded in $B_\epsilon(\tau y_0
+(1-\tau)y_1)$, wherever $u$ is differentiable. This 
contradicts the assumption that $\tau y_0+(1-\tau)y_1\in
\partial_1\Omega$. The theorem follows.
\end{proof}

\begin{remark}
In Theorem \ref{einl thm}, we don't have to assume $\Omega$ 
being connected. If $y_0$ and $y_1$ lie in different 
components of $\Omega$, we find $0<\tau<1$ such that
$\tau y_1+(1-\tau)y_0\in\partial\Omega$. Choose $0<t<1$,
$t\neq\tau$, such that $t y_1+(1-t)y_0\in\Omega$. Assume
that $t>\tau$. There exists $0<\sigma<1$ such that
$$t y_1+(1-t)y_0=\sigma(\tau y_1+(1-\tau)y_0)+(1-\sigma)y_1.$$
Thus
$$u(t y_1+(1-t)y_0)\le\sigma u(\tau y_1+(1-\tau)y_0)
+(1-\sigma)u(y_1)$$
and $u(\tau y_1+(1-\tau)y_0)$ is bounded below. A similar
argument bounds $u$ from below in a neighborhood (relative
to $\Omega$) of $\tau y_1+(1-\tau)y_0$. Following the
lines of the proof of Theorem \ref{main thm}, we obtain
an upper bound on $u$ near $ty_1+(1-t)y_0$ and gradient
bounds. As $\nabla u$ is bounded near $\tau y_1+(1-\tau)y_0$,
we obtain a contradiction. We conclude that $\Omega$ has to be 
connected.
\end{remark}

\begin{remark}
Let $\Omega$ be any convex domain in $\R^n$ and 
$d:=\dist(\cdot,\,\partial\Omega)$. As $d$ is concave in
$\Omega$, $-\sqrt d$ is a convex function in $\Omega$
with unbounded gradient along $\partial\Omega$. 
Thus, according to Theorem \ref{einl thm}, such functions
exist precisely on convex domains.
\end{remark}


\end{document}